\documentclass[12pt,a4]{amsart}
\usepackage[a4paper, left=28mm, right=28mm, top=28mm, bottom=34mm]{geometry}
\usepackage[all]{xy}
\usepackage{enumerate}
\usepackage{amsmath}
\usepackage{amssymb}
\usepackage{amsthm}
\usepackage{mathrsfs}

\usepackage[OT2,T1]{fontenc}
\DeclareSymbolFont{cyrletters}{OT2}{wncyr}{m}{n}
\DeclareMathSymbol{\Sha}{\mathalpha}{cyrletters}{"58}

\newtheorem{thm}{Theorem}[section]
\newtheorem{lem}[thm]{Lemma}
\newtheorem{cor}[thm]{Corollary}
\newtheorem{prop}[thm]{Proposition}

\theoremstyle{definition}

\newtheorem{rem}[thm]{Remark}
\newtheorem{claim}[thm]{Claim}
\newtheorem{defn}[thm]{Definition}

\newtheorem{ex}[thm]{Example}

\numberwithin{equation}{section}

\def\F{{\mathbb F}}

\def\Q{{\mathbb Q}}

\def\Z{{\mathbb Z}}

\def\sO{{\mathscr O}}
\def\sF{{\mathscr F}}

\def\cO{{\mathcal O}}
\def\cA{{\mathcal A}}

\def\res{{\mathrm{res}}}

\def\Abs{\boldsymbol{(\mathrm{Abs})}}
\def\NT{\boldsymbol{(\mathrm{NT})}}
\def\Full{\boldsymbol{(\mathrm{Full})}}

\def\unram{{\mathrm{ur}}}

\def\cont{{\mathrm{cont}}}

\def\nuim{{\nu_{\mathrm{im}}}}
\def\nuimn{{\nu_{\mathrm{im},n}}}

\def\End{\mathop{\mathrm{End}}\nolimits}

\def\Sel{\mathop{\mathrm{Sel}}\nolimits}

\def\cH{\mathop{\mathcal{H}}\nolimits}
\def\wcH{\mathop{\widetilde{\mathcal{H}}}\nolimits}

\def\cF{\mathop{\mathcal{F}}\nolimits}

\def\>'{\mathop{>'}}

\def\cF{\mathop{\mathcal{F}}\nolimits}
\def\cK{\mathop{\mathcal{K}}\nolimits}

\def\Fin{\mathop{\mathrm{Fin}}\nolimits}

\def\cF{\mathop{\mathcal{F}}\nolimits}

\def\Aut{\mathop{\mathrm{Aut}}\nolimits}

\def\Spec{\mathop{\mathrm{Spec}}\nolimits}
\def\Spf{\mathop{\mathrm{Spf}}\nolimits}

\def\Gal{\mathop{\mathrm{Gal}}\nolimits}
\def\Lie{\mathop{\mathrm{Lie}}\nolimits}
\def\coLie{\mathop{\mathrm{coLie}}\nolimits}

\def\Hom{\mathop{\mathrm{Hom}}\nolimits}

\def\corank{\mathop{\mathrm{corank}}\nolimits}

\def\Ker{\mathop{\mathrm{Ker}}\nolimits}
\def\id{\mathop{\mathrm{id}}\nolimits}
\def\GL{\mathop{\mathrm{GL}}\nolimits}
\def\SL{\mathop{\mathrm{SL}}\nolimits}

\def\GSp{\mathop{\mathrm{GSp}}\nolimits}
\def\Sp{\mathop{\mathrm{Sp}}\nolimits}

\def\fsl{\mathop{\mathfrak{sl}}\nolimits}

\def\rank{\mathop{\text{\rm rank}}\nolimits}

\def\divi{\mathop{\mathrm{div}}}
\def\ord{\mathop{\mathrm{ord}}\nolimits}

\begin{document}

\title[class numbers along a Galois representation]{
Asymptotic lower bound of class numbers along a Galois representation}

\author{Tatsuya Ohshita}

\address{Graduate School of Science and Engineering, Ehime University 2--5, Bunkyo-cho, Matsuyama-shi, Ehime 790--8577, Japan}
\email{ohshita.tatsuya.nz@ehime-u.ac.jp}

\date{\today}
\subjclass[2010]{Primary 11R29; Secondary 11G05, 11G10, 11R23. }
\keywords{class number; Galois representation; elliptic curve; 
abelian variety; Selmer group; 
Mordell--Weil group; Iwasawa theory}

\begin{abstract}
Let $T$ be a free $\Z_p$-module of finite rank 
equipped with a continuous $\Z_p$-linear action of 
the absolute Galois group of a number field $K$ 
satisfying certain conditions.   
In this article, by using a Selmer group corresponding to  $T$, 
we give a lower bound of the additive $p$-adic valuation of 
the class number of $K_n$, which is
the Galois extension field of $K$ fixed by 
the stabilizer of $T/p^n T$.
By applying  this result, 
we prove an asymptotic inequality which 
describes an explicit lower bound of 
the class numbers along a tower 
$K(A[p^\infty]) /K$ for a given abelian variety $A$ 
with certain conditions 
in terms of the Mordell--Weil group. 
We also prove another  asymptotic inequality  
for the cases when $A$ is a Hilbert--Blumenthal or 
CM abelian variety.
\end{abstract}

\maketitle

\section{Introduction}\label{secintro}

Commencing with Iwasawa's class number formula 
(\cite{Iw} \S 4.2), 
it is a classical and important problem to study 
the asymptotic behavior of class numbers along 
a tower of number fields. 
Greenberg (\cite{Gr}) and 
Fukuda--Komatsu--Yamagata (\cite{FKY})
studied Iwasawa's $\lambda$-invariant 
of a certain (non-cyclotomic) $\Z_p$-extension of a CM field 
for a prime number $p$:  
by using Mordell--Weil group of 
a CM abelian variety, 
they gave a lower bound of 
the $\lambda$-invariant. 
Sairaiji--Yamauchi (\cite{SY1}, \cite{SY2})
and Hiranouchi (\cite{Hi}) 
studied asymptotic behavior of class numbers along 
a $p$-adic Lie extension $\Q (E[p^\infty])/\Q$ 
generated by coordinates of all $p$-power torsion points 
of an elliptic curve $E$ defined over $\Q$
satisfying certain conditions, and obtained 
results analogous to those in \cite{Gr} and \cite{FKY}.
In this article, by using the terminology of 
Selmer groups, we generalize their results
to the $p$-adic Lie extension of a number field $K$ 
along a $p$-adic representation 
of the absolute Galois group $G_K:=\Gal (\overline{K}/K)$ 
(Theorem \ref{thmmain}). 
As an application of this theory, 
we prove an asymptotic inequality which 
gives a lower bound of 
the class numbers along a tower 
$K(A[p^\infty]) /K$ for a given abelian variety $A$ 
with certain conditions (Corollary \ref{corab}). 
We also prove another  asymptotic inequality  
for the cases when $A$ is a Hilbert--Blumenthal or 
CM abelian variety (Corollary \ref{corCM}).

Let us introduce our notation. 
Fix a prime number $p$, and 
$\ord_p \colon \Q^\times \longrightarrow \Z$ 
the additive $p$-adic valuation normalized by 
$\ord_p (p)=1$. 
Let $K/\Q$ be a finite extension, and
$\Sigma$ a finite set of places of $K$
containing all places above $p$ and 
all infinite places. 
We denote by $K_\Sigma$ the maximal
Galois  extension field of $K$
unramified outside $\Sigma$, and  
put $G_{K,\Sigma}:=\Gal(K_\Sigma/K)$. 
Let $d \in \Z_{>0}$, and 
suppose that a free $\Z_p$-module $T$
of rank $d$ equipped with 
a continuous $\Z_p$-linear $G_{K,\Sigma}$-action 
\(
\rho\colon 
G_{K,\Sigma} \longrightarrow \Aut_{\Z_p} (T)
\simeq \GL_d (\Z_p)
\). 
We put $V:=T \otimes_{\Z_p} \Q_p$, and 
$W:=V/T \simeq T \otimes_{\Z_p} \Q_p/\Z_p$. 
Let $n \in \Z_{\ge 0}$. We denote 
the $\Z_p[G_{K,\Sigma}]$-submodule of $W$
consisting of all $p^n$-torsion elements
by $W[p^n]$. 
We define $K_n:=K(W[p^n])$ to be the maximal subfield of $K_{\Sigma}$
fixed by the kernel of the continuous group homomorphism
\(
\rho_n \colon 
G_{K,\Sigma} \longrightarrow \Aut_{\Z/p^n \Z} (W[p^n])
\)
induced by $\rho$. 
We denote by $h_n$ the class number of $K_n$.
In this article, we study the asymptotic behavior
of the sequence $\{ \ord_p (h_n) \}_{n \ge 0}$ by using 
a Selmer group of $W$. 

Let us introduce some notation related to  Selmer groups in our setting 
briefly. (For details, see \S \ref{ssloccond}.) 
Let 
\(
\cF=\{ H^1_{\cF} (L_w, V) \subseteq H^1 (L_w, V) \}_{L,w}
\) 
be any local condition  on $(V, \Sigma)$
in the sense of Definition \ref{defLC}. 
For instance, we can set $\cF$ to be 
Bloch--Kato's finite local condition $f$.
Let $v \in \Sigma$ be any element, and 
$H^1_{\cF} (K_v, W)$ be the $\Z_p$-submodule of 
$H^1 (K_v, W)$ attached to $\cF$.
Since the Galois cohomology 
$H^1 (K_v, W)$ is a cofinitely generated 
$\Z_p$-module, so is the subquotient 
\[
\cH_v:= H^1_{\cF} (K_v, W)/(H^1_{\cF} (K_v, W) \cap 
H^1_{\unram} (K_v, W)), 
\] 
where 
$H^1_{\unram} (K_v, W)$ denotes the unramified part of 
$H^1 (K_v, W)$.
We denote by the corank of 
the $\Z_p$-module $\cH_v$ by 
$r_v=r_v (T, \cF)$. 
We define the Selmer group of $W$ over $K$
with respect to the local condition $\cF$ by
\[
\Sel_{\cF}(K,W):= \Ker 
\left(
H^1(K_\Sigma /K,W) \longrightarrow 
\prod_{v \in \Sigma} \frac{H^1 (K_v, W)}{H^1_{\cF} (K_v, W)}
\right).
\]
Since $H^1(K_\Sigma /K,W)$ is a cofinitely generated 
$\Z_p$-module, so is $\Sel_{\cF}(K,W)$. 
We define $r_{\Sel}:=r_{\Sel}(T,\cF)$ to be 
the corank of the $\Z_p$-module 
$\Sel_{\cF}(K,W)$. 

For two sequences $\{ a_n \}_{n \ge 0}$ and 
$\{ b_n \}_{n \ge 0}$ of real numbers, 
we write $a_n \succ b_n$ if we have
$\liminf_{n \to \infty} (a_n -b_n)> - \infty$, 
namely if the sequence $\{ a_n -b_n \}_{n \ge 0}$
is bounded below. 
The following theorem is 
the main result of our article.

\begin{thm}\label{thmmain}
Assume that $T$ satisfies the following two conditions. 
\begin{itemize} \setlength{\leftskip}{3mm}
\item[$\Abs$] The representation
$W[p]$ of $G_{K,\Sigma}$ over $\F_p$ is 
absolutely irreducible. 
\item[$\NT$] If $d=1$, then $G_{K,\Sigma}$ acts on $W[p]$ non-trivially. 
\end{itemize}
Then, we have 
\[
\ord_p (h_n) \succ  d \left(r_{\Sel} - 
\sum_{v \in \Sigma}r_v \right) n. 
\] 
\end{thm}

\begin{rem}\label{remstrver}
In this article, we also show a stronger assertion 
than Theorem \ref{thmmain} which describes 
not only asymptotic behavior but also 
a lower bound of each $h_n$ 
in the strict sense. (See  Theorem \ref{thmmainstr}.)  
\end{rem}

Let $A$ be an abelian variety defined over a number field $K$, 
and put $g:=\dim A$. 
For each $N \in \Z_{\ge 0}$, we denote by $A[N]$
the $N$-torsion part  of 
$A(\overline{K})$.  
Put $r_{\Z}(A):=\dim_\Q (A(K)\otimes_{\Z} \Q)$. 
For each $n \in \Z_{\ge 0}$, we denote by $h_{n}(A;p)$
the class number of $K(A[p^n])$, which is the extension field of $K$ 
generated by the coordinates of elements $A[p^n]$. 
By applying Theorem \ref{thmmain}, 
and we obtain the following corollary. 
(See \S \ref{secab}.)

\begin{cor}\label{corab}
Suppose that $A[p]$ is 
an absolutely irreducible representation of 
$G_{K }$ over $\F_p$. 
Then, it holds that  
\[
\ord_p 
(h_n(A;p))
\succ 2g \left(  r_\Z (A)
-g [K : \Q]  \right) n.
\]  
\end{cor}

\begin{rem}
Suppose that $K=\Q$, and $A$ is an elliptic curve 
over $\Q$. 
Then, it follows from Corollary \ref{corab}  that 
\(
\ord_p 
(h_n(A;p))
\succ 2 \left( r_{\Z} (A) 
- 1  \right) n
\). 
This asymptotic inequality coincides with 
that obtained by Sairaiji--Yamauchi
(\cite{SY1}, \cite{SY2}) and 
Hiranouchi (\cite{Hi}). 
Moreover, if $p$ is odd, then  
Theorem \ref{thmmainstr}, which is 
a stronger result than Theorem \ref{thmmain} mentioned in 
Remark \ref{remstrver}, 
implies strict estimates obtained in 
\cite{SY1}, \cite{SY2} and \cite{Hi}. 
For details, see Example \ref{exellstr}.  
\end{rem}

We have another application of Theorem \ref{thmmain} 
for the cases when $A$ is 
a Hilbert--Blumenthal abelian variety or 
a CM abelian variety. 
Suppose that $K$ is a totally real or CM field which is Galois over $K$, 
and denote by $\cO_K$ the ring of integers in $K$. 
We denote by $K^+$ the maximal totally real subfield of $K$, and 
put $g:=[K^+: \Q]$. Let $p$ be a prime number 
which splits completely in $K$, and 
has prime ideal decomposition
\(
p\cO_K = \prod_{\sigma \in \Gal (K/\Q)} \sigma(\pi) \cO_K
\)
for some element $\pi \in \cO_K$. 
Let $A$ be a $g$ dimensional Hilbert--Blumenthal 
(resp.\ CM) abelian variety over $K$  
which has good reduction at every places above $p$, and 
satisfies $\End_K (A) =\cO_K$ if $K$ is totally real 
(resp.\ CM). For each $n \in \Z_{\ge 0}$, 
we denote by $h_{n}(A;\pi)$
the class number of $K(A[\pi^n])$. 
We put $r_{\cO_K}(A):=\dim_K (A(K)\otimes_{\cO_K}K)$.
Then, the following holds.

\begin{cor}\label{corCM}
Suppose that $A[\pi]$ is an absolutely irreducible 
non-trivial representation of $G_K$ over $\F_p$. Then, 
we have 
\[
\ord_p(h_{n}(A;\pi)) \succ \frac{2}{ [K:K^+]}\left(
r_{\cO_K}(A)  -g
\right) n. 
\] 
\end{cor}

\begin{rem}
When $K$ is a CM field,   
the asymptotic inequality in Corollary \ref{corCM}
coincides with 
that obtained by 
Greenberg  \cite{Gr} and 
Fukuda--Komatsu--Yamagata \cite{FKY}.
\end{rem}

The strategy for the proof of our main result, namely
Theorem \ref{thmmain}, is quite similar to 
those established in earlier works \cite{Gr}, \cite{FKY}, 
\cite{SY1}, \cite{SY2} and \cite{Hi}. 
For each $n \in \Z_{\ge 0}$, 
by using elements of the Selmer group,
we construct 
a finite abelian extension $L_n/K_n$ 
which is unramified outside $\Sigma$, 
and whose degree is a power of $p$. 
Then, we compute the degree $[L_n : K_n]$
and the order of inertia subgroups 
at ramified  places.

In \S \ref{secGalcoh}, 
we introduce  notation related to  
Galois cohomology, and prove 
some  preliminary results. 
In \S \ref{secproof}, we prove 
our main theorem, namely Theorem \ref{thmmain}.  
In \S \ref{secab}, we apply Theorem \ref{thmmain} 
to the Galois representations arising from abelian varieties, 
and prove Corollary \ref{corab} and Corollary \ref{corCM}.
We also compare our results with 
earlier works by Fukuda--Komatsu--Yamagata \cite{FKY}, 
Sairaiji--Yamauchi \cite{SY2} and Hiranouchi \cite{Hi}.

\subsubsection*{Notation}
Let $L/F$ be a Galois extension, 
and $M$ a topological abelian group equipped 
with a $\Z$-linear action of $G$.
Then, for each $i \in \Z_{\ge 0}$, 
we denote by 
$H^i(L/F , M):=H^i_{\mathrm{cont}}(\Gal (L/F), M)$
the $i$-th continuous Galois cohomology.
When $L$ is a separable closure of $F$, then
we write $H^i(F , M):=H^i(L/F , M)$. 

Let $F$ be a non-archimedean local field. 
We denote by $F^{\unram}$ the maximal unramified 
extension of $F$.  
For any topological abelian group $M$ equipped with 
continuous $\Z$-linear action of $G_{F}$, we define 
\(
H^1_{\unram}(F,M):= \Ker \left(
H^1(F,M) 
\longrightarrow 
H^1(F^{\unram},M)
\right)
\).

Let $R$ be a commutative ring, and $M$ an $R$-module $M$.  
We denote by $\ell_R (M)$ the length of $M$. 
For each $a \in R$, we denote $M[a]$ 
by the $R$-submodule of $M$ consisting of elements annihilated by $a$. 

\section*{Acknowledgment}
The author would like to thank Takuya Yamauchi 
for giving information on his works with Fumio Sairaiji 
(\cite{SY1}, \cite{SY2}),
and suggesting related problems on abelian varieties. 
This work is motivated by Yamauchi's suggestion.

\section{Preliminaries on Galois cohomology}\label{secGalcoh}

Here, we introduce some notation 
related to Galois cohomology,
and prove preliminary results. 
Let $K$, $\Sigma$ and $T$ be as in \S \ref{secintro}, and  
assume that  $T$ satisfies the conditions $\Abs$ and $\NT$
in Theorem \ref{thmmain}.
We denote by $\Fin (K;\Sigma )$
be the set of all intermediate fields $L$
of $K_\Sigma /K$ which are finite over $K$. 
For each $L \in \Fin (K;\Sigma )$, we denote by 
$P_L$ the set of all places of $L$, and 
by $\Sigma_L$
the subset of $P_L$ consisting of places above an element of $\Sigma$.

\subsection{Local conditions}\label{ssloccond}
In this subsection, 
let us define the notion of local conditions 
and Selmer groups in our article.

\begin{defn}\label{defLC}
Recall that we put $V:= T \otimes_{\Z_p} \Q_p$, 
and $W:=V/T$.
\begin{enumerate}[{\rm (i)}]\setlength{\parskip}{1mm}
\item A collection 
\(
\cF :=
\left\{ H^1_{\cF}(L_w,V) \subseteq H^1(L_w , V)  
\mathrel{\vert} L \in \Fin (K; \Sigma) ,\, w \in \Sigma_L 
\right\}
\)
of $\Q_p$-subspaces is called {\em a local condition} on $(V, \Sigma)$ 
if  the following $(*)$ is satisfied.  
\begin{enumerate}\setlength{\parskip}{1mm}
\item[$(*)$]  {\em Let $\iota \colon L_1 \hookrightarrow  L_2 $ 
be an embedding of fields 
belonging to $\Fin_K (L; \Sigma)$ over $K$. 
Then, for any $w_1\in P_{L_1}$ and $w_2 \in P_{L_2}$ 
satisfying $\iota^{-1} w_2 =w_1$, the image of 
$ H^1_{\cF}(L_{1,w_1},V)$ via the map 
$ H^1 (L_{1,w_1},V) \longrightarrow  H^1 (L_{2,w_2},V)$
induced by $\iota$ is contained in $H^1_{\cF} (L_{w_2},V)$.}
\end{enumerate}
\item Let $L \in \Fin (K; \Sigma)$ and $w \in P_L$. 
Then, we define $ H^1_{\cF}(L_w,W)$ to be the image of 
$ H^1_{\cF}(L_w,V)$ via the natural map 
$ H^1(L_w,V) \longrightarrow  H^1(L_w,W)$. 
For any $n \in \Z_{\ge 0}$, 
we define  $H^1_{\cF}(L_w,W[p^n])$ to be the inverse image of 
$ H^1_{\cF}(L_w,W)$ via the natural map 
$ H^1(L_w,W[p^n]) \longrightarrow  H^1(L_w,W)$. 
\item Let $L \in \Fin (K; \Sigma)$, and 
$n \in \Z_{\ge 0} \cup \{ \infty \}$.  
Then, we define
\[
\Sel_{\cF}(L,W[p^n]):= 
\Ker 
\left(
H^1(K_\Sigma /L,W[p^n]) \longrightarrow 
\prod_{w \in \Sigma_L} \frac{H^1 (L_w, W[p^n])}{H^1_{\cF} (L_w, W[p^n])}
\right).
\]
\end{enumerate}
\end{defn}

\begin{rem}
Let $\cF$ be a local condition on $(V,\Sigma)$. 
Then, by definition, for any 
$L \in \Fin (K; \Sigma)$ and $w \in P_L$, 
the $\Z_p$-module $H^1_{\cF}(L_w, W)$ is divisible. 
\end{rem}

\begin{rem}
Let $L \in \Fin (K; \Sigma)$ any element, 
and $w \in \Sigma_L$ an infinite place. 
Then, we note that $H^1 (L_w, V)=0$.
Thus for any local condition $\cF$ on $(V, \Sigma)$, 
it clearly holds that $H^1_{\cF} (L_w, V)=0$. 
We also note that $H^1 (L_w, W)$ is annihilated by $2$. 
In particular $H^1 (L_w, W)$ never has 
a non-trivial divisible $\Z_p$-submodule, 
So, the corank of $H^1 (L_w, W)$ is zero.
When we treat a local condition, 
we may not care infinite places.  
\end{rem}

\begin{ex}[\cite{BK} \S 3]

For $L \in \Fin (K; \Sigma)$ and finite place $w \in P_L$, we put
\[
H^1_f (L_w, V):= 
\begin{cases}
H^1_{\unram} (L_w, V) &  (\text{if $w \nmid p$}), \\
\Ker \left(  
H^1 (L_w, V) \longrightarrow 
H^1 (L_w, V\otimes_{\Q_p} B_{\mathrm{crys}}) 
\right)
&  (\text{if $w \mid p$}), \\
0 & (\text{if $w \mid \infty$})
\end{cases}
\]
where $B_{\mathrm{crys}}$ is Fontaine's $p$-adic period ring 
introduced in \cite{Fo} and \cite{FM}. 
Then, we can easily verify that the collection 
\(
\{ 
H^1_f (L'_{w'}, V) \mathrel{\vert}
L' \in \Fin (K; \Sigma), \ w' \in \Sigma_{L'}
\}
\)
forms a local condition on $(V,\Sigma)$.
We call this collection 
{\em Bloch--Kato's finite local condition}.
\end{ex}

\subsection{Global cohomology}

In this subsection, we introduce some 
preliminaries on global Galois cohomology.
We put $K_\infty:=\bigcup_{n \ge 0} K_n$. 
For each $m,n \in \Z_{\ge 0} \cup \{ \infty \}$ with $n \ge m$,
we put
$G_{n,m} :=\Gal(K_n/K_m)$, and 
$G_{m} :=\Gal(K_\Sigma /K_m)$.

First, we control the kernel of 
the restriction map 
\[
\res_{n,W} \colon H^1 (K,W[p^n]) \longrightarrow 
H^1 (K_n,W[p^n])
\]
for every $n \in \Z_{\ge 0}$. 
In order to do it, the following fact 
and the irreducibility of $V$, 
which follows from $\Abs$ for $T$
become a key.

\begin{thm}[\cite{Ru} Theorem C.1.1 in Appendix C]\label{thmvanish}
Let $V'$ be a finite dimensional $\Q_p$-vector space, 
and suppose that a compact subgroup $H$ of 
$\GL (V'):=\Aut_\Q (V')$ acts irreducibly on $V'$
via the standard action. 
Then, we have $H^1 (H, V')=0$. 
\end{thm}

By using Theorem \ref{thmvanish}, 
let us prove the following lemma. 

\begin{lem}\label{lemGlcohfin}
There exist a positive integer $\nuim$ 
such that for any $n \in \Z_{\ge 0}$, 
we have 
$\ell_{\Z_p}(\Ker \res_{n,W}) \le \nuim$.
\end{lem}

\begin{proof}
Let us show that 
$\# H^1 (K_\infty /K, W)< \infty$.
By  $\Abs$ and Theorem \ref{thmvanish}, we have
$H^1 (K_\infty /K, V)=0$.
So, we obtain an injection 
$H^1 (K_\infty /K, W) \hookrightarrow 
H^2 (K_\infty /K, T)$. 
Note that 
$G_{\infty,0}=\Gal(K_\infty / K)$ 
is a compact $p$-adic analytic group 
since $G_{\infty,0}$ can be regarded as 
a compact subgroup of $\GL (V)$.
So, by Lazard's theorem (for instance, see \cite{DDMS} 8.1 Theorem),
it holds that  
$G_{\infty,0}$ is topologically finitely generated. 
This implies that the order of $H^2 (K_\infty/K, W[p])$ is finite.
Thus $H^2 (K_\infty /K, T)$ is finitely generated over $\Z_p$. 
(See \cite{Ta} Corollary of (2.1) Proposition.)
Since $H^1 (K_\infty /K, W)$ is a torsion $\Z_p$-module, 
we deduce that the order of 
$H^1 (K_\infty /K, W)$
is finite.

Let $\nu$ be the length of 
the $\Z_p$-module $H^1 (K_\infty /K, W)$. 
Take any $n \in \Z_{\ge 0}$. 
By the assumptions $\Abs$ and $\NT$, 
the natural map 
$H^1 (K_n /K, W[p^n]) \longrightarrow H^1 (K_n /K, W)$ 
is injective. 
So, by the inflation-restriction exact sequence, 
we obtain an injection
\(
\Ker \res_{n,W}
\hookrightarrow 
H^1 (K_\infty /K, W)
\). 
Hence $\nuim :=\nu$ satisfies the desired properties. 
\end{proof}

\begin{rem}\label{remnuimvanish}
Note that if the image 
$\rho_1 (\Gal (K_1/K)) \subseteq \GL_d (\F_p)$
contains a non-trivial scalar matrix, then 
we can take $\nuim =0$. Indeed, in such cases, 
similarly  to \cite{LW} \S 2 Lemma 3, 
we can show that  
$H^1 (K_\infty /K, W)=0$.
\end{rem}

\begin{ex}\label{exmnuimvanish}
Suppose $d= \dim_{\Q_p} V=2$, and 
the following assumption $\Full$. 
\begin{itemize}\setlength{\leftskip}{3mm}
\item[$\Full$] If $p$ is odd, then  the map
$\rho_1 \colon G_{K,\Sigma} \longrightarrow 
\Aut_{\F_p}(W[p]) \simeq \GL_2 (\F_p)$ 
is surjective. 
If $p=2$, then $\rho \colon G_{K,\Sigma} \longrightarrow 
\Aut_{\Z_2}(T) \simeq \GL_2 (\Z_2)$ is surjective. 
\end{itemize}
Then, we may take $\nuim=0$ 
because of the following Claim \ref{claimvanish}. 
Note that in 
\cite{SY1}, \cite{SY2} and \cite{Hi}, 
Sairaiji, Yamauchi and Hiranouchi assumed 
the hypothesis $\Full$.
So, the constant $\nuim$ does not appear explicitly in these works.
\end{ex}

\begin{claim}\label{claimvanish}
If $d=2$, and if $T$ satisfies $\Full$, then we have 
$H^1 (K_\infty /K, W)=0$.
\end{claim}

\begin{proof}[Proof of Claim \ref{claimvanish}]
If $p$ is odd, then the claim follows from 
similar arguments to those in the proof of 
\cite{LW} \S 2 Lemma 3. 
So, we may suppose that $p=2$. 
It suffices to show that 
$H^0(K, H^1 (K_{n+1} /K_n, W[2]))=0$ 
for each $n \in \Z_{\ge 0}$.
Note that we have $G_{1,0} \simeq \GL_2 (\F_2) \simeq \mathfrak{S}_3$.
Let $A$ be a unique normal subgroup of $G_{1,0}$ of order $3$. 
Then, we have $W[2]^A=0$, and $H^1 (A, W[2])=0$. 
So, we have  $H^1 (K_{1} /K, W[2])=0$. 

Take any $n\ge 1$, and let us show that  
\begin{equation}\label{eqKn+1Knvanish}
H^0(K, H^1 (K_{n+1} /K_n, W[2]))=\Hom(G_{n+1,n}, 
W[2])^{G_{n+1,0}}=0.  
\end{equation}
The map $\rho_{n+1} \colon G_{n+1,0}
\longrightarrow \GL_2(\Z/2^{n+1} \Z)$ 
induces an isomorphism from $G_{n+1, n}$  to 
a subgroup of 
\(
(1+2^nM_2(\Z_2))/((1+2^{n+1}M_2(\Z_2))) 
\simeq M_2(\F_2)
\)
preserving the conjugate action of 
$G_{n+1,0} \simeq \rho_{n+1}(G_{n+1, 0})$, 
which factors through $G_{1,0} \simeq \GL_2 (\F_p)$.
By the assumption $\Full$, we have 
$\Gal (K_{n+1}/K_n) \simeq M_2 (\F_2)$. 
Note that the $\F_2 [\GL (\F_2)]$-submodules of 
$M_2 (\F_2)$ are $0 \subseteq \F_2 \subseteq 
\fsl (\F_2) \subseteq M_2 (\F_2)$. 
So, we deduce that $M_2 (\F_2)$ never has a quotient 
isomorphic to $\F_2^2$. 
Hence the equality (\ref{eqKn+1Knvanish})  holds. 
\end{proof}

Next, by using Galois cohomology classes 
contained in the Selmer group,
we shall construct certain number fields. 
Let $n \in \Z_{\ge 0}$ be any element. 
Clearly, we have a natural isomorphism
\(
H^1(K_\Sigma/K_n,W[p^n])\simeq 
\Hom_{\cont}(G_n, W[p^n])  
\),  
where $\Hom_{\cont}(G_n, W[p^n])$
denotes the group consisting of 
continuous homomorphisms from 
$G_n$  to $ W[p^n]$.
Since $G_n$ is a normal subgroup 
of $G_{K,\Sigma}=G_0$, 
we can define a left action  
\[
G_{n,0} \times 
\Hom_{\cont}(G_n, W[p^n])
 \longrightarrow 
\Hom_{\cont}(G_n, W[p^n]);\ 
(\sigma, f)  \longmapsto \sigma * f
\]
of $G_{n,0}$ on $\Hom_{\cont}(G_n, W[p^n])$ 
by 
\(
(\sigma * f)(x):=\sigma 
(f(\widetilde{\sigma}^{-1} x \widetilde{\sigma}))
\) 
for each $x \in G_n$, where 
$\widetilde{\sigma} \in G_0$ is a lift of $\sigma$.
Note that he definition of $\sigma * f$ is independent of  
the choice of $\widetilde{\sigma}$. 
The following lemma is a key of the proof of 
Theorem \ref{thmmain}.

\begin{lem}\label{lemlengthdeg}
Take any $n \in \Z_{\ge 0}$. 
Let $M$ be a $\Z_p$-submodule of 
\[
\cH_n:= 
H^0(K, H^1(K_\Sigma/K_n,W[p^n]))=\Hom_{\cont}(G_n, 
W[p^n])^{G_{n,0}}. 
\]
We define $K_n (M)$ to be the maximal subfield of $K_\Sigma$
fixed by $\bigcap_{h \in M} \Ker h$. 
Then $K_n (M) /K_n$ is Galois,  
and 
$[K_n (M): K_n]= p^{d \ell_{\Z_p} (M)}$. 
Moreover, the evaluation map 
\(
e_{M} \colon M  \longrightarrow \Hom_{\Z_p [G_{n,0}]} 
\left(
\Gal (K_n (M)/K_n), 
W[p^n] \right) 
\)
is  an isomorphism of $\Z_p$-modules.
\end{lem}

\begin{proof}
By definition, the extension $K_n (M) /K_n$ is 
clearly Galois, and $e_M$ is a well-defined injective homomorphism.  
Let us show the rest of the assertion of Lemma \ref{lemlengthdeg}
by induction on $\ell_{\Z_p} (M)$.

When $\ell_{\Z_p} (M)=0$, 
the assertion of Lemma \ref{lemlengthdeg} is clear.

Let $\ell$ be a positive integer, and 
suppose that the assertion of Lemma \ref{lemlengthdeg} holds 
for any $\Z_p$-submodule $M'$ of 
$\cH_n$ satisfying $\ell_{\Z_p} (M') < \ell$. 
Let $M$ be any 
$\Z_p$-submodule of $\cH_n$ 
satisfying $\ell_{\Z_p} (M) = \ell$. 
Take a $\Z_p$-submodule $M_0$ of $M$
such that $\ell_{\Z_p} (M_0) = \ell-1$. 
By definition, we have $K_n (M_0) \subseteq K_n(M)$.
Since $e_{M_0}$ is an isomorphism 
by the hypothesis of induction,
and since $e_M$ is an injection, 
we deduce that 
\begin{equation}\label{eqdeg>1}
[K_n (M): K_n (M_0)]>1. 
\end{equation}
For each $\Z_p$-submodule $N$ of $M$, we put 
$\mathfrak{K}(N):= \bigcap_{h' \in N} \Ker h'
= \Gal (K_\Sigma/K_n({N}))$. 
Take an element $f \in M$ not contained in $M_0$. 
Then, we have 
$\mathfrak{K}(M)=\Ker f 
\cap \mathfrak{K}(M_0)$. 
Note that the abelian group 
$\Gal (K_n (M)/ K_n (M_0))=
\mathfrak{K}(M_0)/\mathfrak{K}(M)$
is annihilated by $p$.  
(Indeed, if there exists an element 
$\Gal (K_n (M)/ K_n (M_0))$
which is not annihilated by $p$, then 
we obtain a sequence $M_0 \subset pM+M_0 \subset M$,  
which contradicts the fact that $M/M_0$
is a simple $\Z_p$-module.) 
So, the map $f$ induces an injective  
$\F_p [G_{n,0}]$-linear map from 
$\Gal (K_n (M)/ K_n (M_0))=\mathfrak{K}(M_0)/
(\Ker f \cap \mathfrak{K}(M_0))$ into 
$(W[p^n])[p]=W[p]$.
By the inequality (\ref{eqdeg>1}) and the assumption $\Abs$,  
we deduce that 
\begin{equation}\label{eqGalsimple}
\Gal (K_n (M)/ K_n (M_0)) \simeq W[p].
\end{equation}
Since we have $[K_n (M_0) : K_n]=p^{d(\ell-1)}$ 
by the induction hypothesis, 
we obtain 
\[
[K_n (M): K_n]= p^{d(\ell-1)}\cdot \# W[p] 
=p^{d \ell}  =
p^{d \ell_{\Z_p} (M)}.
\]
In order to complete the proof of Lemma \ref{lemlengthdeg}, 
it suffices to prove that the map $e_M$ is surjective. 
For each $\Z_p$-submodule $N$ of $M$, we put 
\[
X(N):=
\Hom_{\Z_p [G_n]} 
\left(
\Gal (K_n (M)/K_n), 
W[p^n] \right).
\]
Since $e_M$ is injective, and since 
$\ell_{\Z_p} (M)=\ell$, 
it suffices to show that 
$\ell_{\Z_p} (X(M)) \le \ell$. 
By the induction hypothesis, 
we have $\ell_{\Z_p} (X(M_0))=\ell-1$.
Put 
\[
X(M;M_0):=
\Hom_{\Z_p [G_{n,0}]} 
\left(
\Gal (K_n (M)/K_n(M_0)), 
W[p^n] \right).
\]
Since we have an exact sequence 
\(
0  \longrightarrow 
X(M;M_0) \longrightarrow 
X(M) \longrightarrow 
X(M_0)
\),
it suffices to show that 
the $\Z_p$-module $X(M;M_0)$ is simple. 
By  (\ref{eqGalsimple}), we obtain 
\[
X(M;M_0) \simeq \Hom_{\Z_p [G_{n,0}]} (W[p],W[p^n]) 
\simeq \End_{\F_p [G_{n,0}]} (W[p]) . 
\]
Since the representation $W[p]$ of $G_n$
is absolutely irreducible over $\F_p$ by the assumption $\Abs$, 
we obtain $\End_{\F_p [G_{n,0}]} (W[p]) =\F_p$.
Hence the $\Z_p$-module $X(M;M_0)$ is simple.
This completes the proof of Lemma \ref{lemlengthdeg}. 
\end{proof}

\section{Proof of Theorem \ref{thmmain}}\label{secproof}

In this section, we prove Theorem \ref{thmmain}. 
Let us fix our notation. 
Again, let $K$, $\Sigma$ and $T$ be as in \S \ref{secintro}. 
Assume that  $T$ satisfies the condition $\Abs$ and $\NT$.
Take any local condition $\cF$ on $(V,\Sigma)$. 
Recall that we put $r_{\Sel}=r_{\Sel} (T, \cF)
= \corank_{\Z_p} \Sel_{\cF} (K, W)$. 
For each $n \in \Z_{\ge 0}$, we denote by $M_n$ the image of 
\(
\Sel_{\cF} (K, W[p^n]) \longrightarrow 
\Sel_{\cF} (K_n, W[p^n]),
\) 
and we put $L_n:=K_n(M_n)$ in the sense of 
Lemma \ref{lemlengthdeg}.

\begin{prop}\label{propglob}
Let $\nuim $ be as in Lemma \ref{lemGlcohfin}. 
For any w$n \in \Z_{\ge 0}$, 
the extension $L_n/K_n$ is unramified outside $\Sigma_{K_n}$, 
and 
\(
\ord_p 
[L_n : K_n] \ge d (n r_{\Sel}  - \nuim)
\). 
\end{prop}

\begin{proof}
Since every $f \in M_n$ is a map defined on $G_{K,\Sigma}$, 
we deduce that $L_n$ is 
unramified outside $\Sigma$. 
By Lemma \ref{lemGlcohfin}, the length of 
the $\Z_p$-module $M_n$
is at least $n r_{\Sel}  - \nuim$.
So, Lemma \ref{lemlengthdeg} implies that 
$\ord_p [L_n : K_n] \ge d (n r_{\Sel}  - \nuim)$. 
\end{proof}

Let $v \in \Sigma$ be any finite place.
We denote by $P_{n,v}$ 
the set of all places of $K_n$
above $v$. 
For each $w \in P_{n,v}$, 
we denote by $I_{w}(L_n/K_n)$ 
the inertia subgroup 
of $\Gal (L_n/K_n)$ at $w$. 
We define  $I_{n,v}$ to 
be the subgroup of $\Gal (L_n/K_n)$ 
generated by $\bigcup_{w \in P_{n,v}} I_w (L_n /K_n)$. 
Recall that $r_{v}=r_{v} (T, \cF)$ denotes 
the corank of the $\Z_p$-module
\(
\cH_v:= H^1_{\cF} (K_v, W)/(H^1_{\cF} (K_v, W) \cap 
H^1_{\unram} (K_v, W))
\). 
Let $\cK /K_v$ be an algebraic extension.
Put $W(\cK):=H^0 (\cK,W)$. 
We define  
$W(\cK)_{\divi}$ to be the maximal divisible 
$\Z_p$-submodule of $W (\cK)$. 
For each $n \in \Z_{\ge 0} \cup \{ \infty \}$, we define 
\[
\nu_{v,n} := \ell_{\Z_p} 
\left(H^0 (K_v, W(K_v^{\unram}) / W(K_v^{\unram})_{\divi})
\otimes_{\Z} \Z/p^n\Z\right).
\]
Note that $\{ \nu_{v, n} \}_{n \ge 0}$ is a bounded increasing sequence, 
and for any sufficiently large $m$, we have $\nu_{v, m} =\nu_{v,\infty}$.
Let us prove the following proposition.

\begin{prop}\label{proploc}
We have $\ell_{\Z_p}(I_{n,v}) \le 
d(r_v n + \nu_{v, n })$
for any $v \in \Sigma$ 
and  $n \in \Z_{\ge 0}$.
\end{prop}

In order to show that Proposition \ref{proploc}, 
we need the following lemma. 

\begin{lem}\label{lemunrerro}
Let $n \in \Z_{\ge 0}$ be any element, 
and define 
\[
\wcH_{v,n}:= H^1_{\cF} (K_v, W[p^n])/(H^1_{\cF} (K_v, W[p^n]) 
\cap H^1_{\unram} (K_v, W[p^n])).
\]
Then, we have $\ell_{\Z_p} (\wcH_{v,n}) \le 
r_{v}n+ \nu_{v,n}$.
\end{lem}

\begin{proof}[Proof of Lemma \ref{lemunrerro}]
Let $\cK /K_v$ be an algebraic extension, and 
\[
\iota_{\cK ,n} \colon H^1 (\cK ,W[p^n]) 
\longrightarrow H^1 (\cK ,W)
\]
the natural map. 
By the short exact sequence 
$0 \longrightarrow W[p^n] \longrightarrow 
W \longrightarrow W \longrightarrow 0$, 
we obtain isomorphism
\(
\Ker \iota_{\cK, n} \simeq W(\cK) \otimes_{\Z} \Z/p^N\Z
\simeq (W(\cK) / W(\cK)_{\divi})\otimes_{\Z} \Z/p^N\Z.
\) 

We put 
$\widetilde{Y}_n:=H^1_{\cF} (K_v, W[p^n]) 
\cap H^1_{\unram} (K_v, W[p^n])$, and 
$Y:=H^1_{\cF} (K_v, W) 
\cap H^1_{\unram} (K_v, W)$.
By definition, it clearly holds that 
$\widetilde{Y}_{n}\subseteq \iota_{K_v,n}^{-1} (Y)$, 
and 
$H^1_{\cF} (K_v, W[p^n])/\widetilde{Y}_{n} = \wcH_{v,n}$.
Moreover, we have an  injection 
\(
H^1_{\cF} (K_v, W[p^n])/\iota_{K_v,n}^{-1} (Y) 
\hookrightarrow \cH_{v}[p^n]
\). 
So, we obtain 
\[
\ell_{\Z_p}(\wcH_{v,n}) \le \ell_{\Z_p}(\cH_{v}[p^n])+ 
\ell_{\Z_p}(\iota_{K_v,n}^{-1} (Y)/\widetilde{Y}_{n} ).
\] 
On the one hand, 
since $H^1_{\cF} (K_v, W)$ is a divisible $\Z_p$-module by definition, 
so is the quotient module $\cH_{v}$. 
This implies that $\ell_{\Z_p}(\cH_{v}[p^n])= r_v n$. 
On the other hand, the restriction map 
$H^1(K_v, W[p^n]) \longrightarrow H^0(K_v, H^1(K^{\unram}_v, W[p^n]))$
induces an injection
\[
\iota_{K_v,n}^{-1} (Y)/\widetilde{Y}_{n} 
\hookrightarrow H^0(K_v, \Ker \iota_{K^{\unram}_v,n})
\simeq 
H^0(K_v, W(K^{\unram}_v) 
/ W(K^{\unram}_v)_{\divi})
\otimes_{\Z} \Z/p^n\Z.
\]
So, we obtain $\ell_{\Z_p}(\iota_{K_v,n}^{-1} (Y)/\widetilde{Y}_{n} )
\le \nu_{v,n}$. 
Hence $\ell_{\Z_p}(\wcH_{v,n}) \le r_{v}n+ \nu_{v,n}$.
\end{proof}

\begin{proof}[Proof of Proposition \ref{proploc}]
Take any $v \in \Sigma$ and 
$n \in \Z_{\ge 0}$. 
Let $w \in P_{n,v}$. We define
\[
\res_{I,w} \colon \Hom_{\Z_p} (\Gal(L_n/K_n),W[p^n]) 
\longrightarrow 
\Hom_{\Z_p} (I_{w}(L_n/K_n),W[p^n])
\]
to be the restriction maps, and put 
$M^{\unram}_{n,w}:= \Ker 
(\res_{I,w} \circ \res_{D,w} \vert_{M_n}) 
\subseteq M_n$.
By definition, the extension 
$K_n(M^{\unram}_{n,w})/K_n$
is unramified at $w$.
So, we obtain 
\begin{equation}\label{eqIwMur}
I_w (L_n / K_n) \subseteq \Gal (L_n/ K_n(M^{\unram}_{n,w})).
\end{equation}
We fix an element $w_0 \in P_{n,v}$. 
Let $\sigma \in G_n$ be any element. 
Then, the diagram
\[
\xymatrix{
M_n \ar[rr]^{\sigma *(-) =\id_{M_n}} \ar[d]_{\res_{I,\sigma^{-1} w_0}} & &
M_n \ar[d]^{\res_{I, w_0}} \\
\Hom_{\Z_p} (I_{\sigma^{-1}w_0}(L_n/K_n),W[p^n]) 
\ar[rr]^{\sigma *(-)}_{\simeq}  & &
\Hom_{\Z_p} (I_{w_0}(L_n/K_n),W[p^n])
}
\]
commutes.
So, for each $w \in P_{n,v}$, 
we have $M^{\unram}_{n,w}=M^{\unram}_{n,w_0}$. 
Hence by (\ref{eqIwMur}), 
we obtain 
\(
I_{n,v} \subseteq \Gal (L_n/ K_n(M^{\unram}_{n,w_0})) 
\). 
In order to prove Proposition \ref{proploc}, 
it suffices to show that 
\begin{equation}\label{eqineqMnMnur}
\ell_{\Z_p} \left(
\Gal (L_n/ K_n(M^{\unram}_{n,w_0})) \right)
\le d(r_v n + \nu_{v , n}).
\end{equation}
By Lemma \ref{lemlengthdeg}, we have  
\(
\ell_{\Z_p}(\Gal (L_n/ K_n(M^{\unram}_{n,w_0})))
=d\ell_{\Z_p}(M_n/ M^{\unram}_{n,w_0})
\).  
Since the natural surjection  
\(
\xymatrix{
\Sel_{\cF} (K,W[p^n]) \ar@{->>}[r] & M_n/
M^{\unram}_{n,w_0}
\simeq (\res_{I,w} \circ \res_{D,w})(M_n)
}
\)
factors through the $\Z_p$-module $\wcH_{v,n}$
in Lemma \ref{lemunrerro}, 
we obtain the inequality (\ref{eqineqMnMnur}).
\end{proof}

\begin{proof}[Proof of Theorem \ref{thmmain}]
Take any $n \in \Z_{\ge 0}$.
Let $I$ be the subgroup of $\Gal (L_n/K_n)$
generated by $\bigcup_{v \in \Sigma} I_{n,v}$. 
Then, the extension $L_n^I/K_n$ is 
unramified at every finite place, 
and the degree $[L^I_n:K_n]$ is a power of $p$.
So, by the global class field theory, 
we have 
\[
\ord_p (h_n)  \ge \ord_p [L^I_n: K_n] 
 \ge \ord_p [L_n:K_n] - \sum_{v \in \Sigma} \ord_p (\#I_{n,v}).
\]
For each $n \in \Z_{\ge 0}$, 
we put $\nuimn:= \ell_{\Z_p} (H^1(K_n/K, W[p^n]))$. 
Then, by Lemma \ref{lemGlcohfin}, 
the sequence $\{ \nuimn \}_{n \ge 0}$ is bounded.
By Proposition \ref{propglob} and
Proposition \ref{proploc}, we obtain
\begin{equation}\label{eqprec}
\ord_p (h_n) \ge d \left(r_{\Sel} \cdot n- \nuimn \right) 
- d \sum_{v \in \Sigma} \left( r_v  n + \nu_{v, n} \right)
\succ d \left(r_{\Sel} - 
\sum_{v \in \Sigma}r_v \right) n. 
\end{equation}
This completes the proof of Theorem \ref{thmmain}. 
\end{proof}

By the above arguments, 
in particular by the inequality (\ref{eqprec}), 
we have also obtained a bit stronger result than 
Theorem \ref{thmmain}.

\begin{thm}\label{thmmainstr}
For any $n \in \Z_{\ge 0}$, we have 
\[
\ord_p (h_n) \ge d \left(r_{\Sel} - 
\sum_{v \in \Sigma}r_v \right) n
- d 
\nuimn - d\sum_{v \in \Sigma}\nu_{v, n}.
\]
\end{thm}

\section{Application to abelian varieties}\label{secab}

In this section, we apply 
Theorem \ref{thmmain}  
to the extension defined by  an abelian variety $A$.
In \S \ref{subsecnoncm}, 
we prove Corollary \ref{corab}.  
Moreover, from the view point of Theorem \ref{thmmainstr}, 
we compare our results (of stronger form) with 
earlier results 
in the cases when $A$ is an elliptic curve. 
In \S \ref{seccm}, we study
the cases when  $A$ is a Hilbert--Blumenthal or CM abelian variety, 
and prove Corollary \ref{corCM}.

\subsection{General cases}\label{subsecnoncm}

Let $A$ be an abelian variety 
over a number field $K$, 
and fix a  prime number $p$ such that $A[p]$ becomes 
an absolutely irreducible representation of 
the absolute Galois group $G_K$ of $K$ over $\F_p$. 
We denote the dimension of $A$ over $K$ by $g$.
Let $T_pA$ be the $p$-adic Tate module of $A$, 
namely $T_pA :=\varprojlim_n A[p^n]$, 
and put $V_p A := T_p A \otimes_{\Z_p} \Q_p$. 
Note that $T_p A$ is a free $\Z_p$-module of rank $2g$. 
Let $\Sigma(A)$ be the subset of $P_K$ consisting of 
all places dividing $p \infty$ and all places 
where $A$ has bad reduction.
Then, the natural action $\rho_A^{(p)}$ of $G_K$ 
is unramified outside $\Sigma(A)$.

Let $L/K$ be any finite extension, and 
$w \in P_K$ any finite place above a prime number $\ell$.
With the aid of the implicit function theorem (for instance \cite{Se} 
PART II Chapter III \S 10.2 Theorem) and the Jacobian criterion 
(for insatnce \cite{Li} Chapter 4 Theorem 2.19), 
the projectivity and smoothness of $A$ implies that 
$A(L_w)$ is a $g$-dimensional compact abelian analytic group over $L_w$. 
So, we have 
\begin{equation}\label{eqstrthm}
A(L_w) \simeq \Z_\ell^{g[L_w: \Q_\ell]} \oplus (\text{a finite abelian group})
\end{equation}
(See Corollary 4 of Theorem 1 in  \cite{Se} PART II Chapter V \S 7.)
Related to this fact, 
the following is known.

\begin{prop}[\cite{BK} Example 3.11]\label{corclvsfin}
For any finite extension field $L$ of $K$, and 
any finite place $w \in P_L$, we have a natural isomorphism
$H^1_f (L_w, A[p^\infty]) \simeq A(L_w) \otimes_{\Z} \Q_p/\Z_p $. 
If  $w$ lies above $p$, then 
the corank of the $\Z_p$-module $H^1_{f}(L_w,A[p^\infty])$ 
is equal to $g [L_w:\Q_p]$. 
If  $w$ does not lie above $p$, then $H^1_{f}(L_w,A[p^\infty])=0$. 
\end{prop}

\begin{proof}[Proof of Corollary \ref{corab}]
We define  the Tate-Shafarevich group $\Sha (A/K) $ 
to be the kernel of  
\(
H^1(K,A(\overline{K})) \longrightarrow 
\prod_{v \in P_K} H^1 (K_v, A(\overline{K}_v))
\). 
Then, we have a short exact sequence 
\begin{equation}\label{eqexseqSel}
0 \longrightarrow A(K) \otimes_\Z  \Q_p/\Z_p 
\longrightarrow 
\Sel_f(K,A[p^\infty]) \longrightarrow 
\Sha (A/K)[p^\infty] \longrightarrow 0.
\end{equation}
(See \S C.4 in \cite{HS} Appendix C. Note 
that by Proposition \ref{corclvsfin}, our 
$\Sel_f(K,A[p^\infty])$ is naturally isomorphic to 
$\varinjlim_n \Sel^{([p^n])}(A/K)$ in the sense of 
\cite{HS} Appendix C, where $[p^n]\colon A \longrightarrow A$ 
denotes the multiplication-by-$p^n$ isogeny.)
So, we obtain
\[
\corank_{\Z_p} \Sel_f(K,A[p^\infty]) 
\ge r_\Z (A):= \rank_{\Z} A(K). 
\]
Hence by Theorem \ref{thmmain} for 
$(T,\Sigma, \cF)=(T_pA, \Sigma (A), f)$ and 
Proposition \ref{corclvsfin}, 
we obtain  
\[
\ord_p 
(h_n(A;p))
 \succ 
2g \left( r_{\Z}(A)  
- g \sum_{v \mid p}[K_v : \Q_p]  \right) n 
=2g \left(r_{\Z} (A) - g [K : \Q]  \right) n.
\]
This completes the proof of  Corollary \ref{corab}. 
\end{proof}

\begin{rem}\label{remnuimforab}
Here, we give remarks on 
the image of the modulo $p$ representation  
\(
\rho_{A, 1}^{(p)} \colon 
\Gal (K(A[p])/K) 
\longrightarrow \Aut (A[p]) 
\). 
Let $A$ be a principally polarized abelian variety 
of dimension $g$ defined over $K$, and 
take an odd prime number $p$. 
Then, the image of $\rho_1$ can be regarded as a subgroup of 
$\GSp_{2g} (\F_p)$. 
Clearly if $\mathrm{Im}\, \rho_{A,1}^{(p)}$ contains 
$\Sp_{2g} (\F_p)$, 
then $T_p A$ 
satisfies the conditions $\Abs$ and $\NT$. 
Since $p$ is odd, the non-trivial scalar $-1$ 
is contained in $\Sp_{2g} (\F_p)$. 
So, as noted in Remark \ref{remnuimvanish}, if 
$\mathrm{Im}\, \rho_{A,1}^{(p)}$ contains 
$\Sp_{2g} (\F_p)$, then we can take $\nuim=0$, 
where $\nuim$ denotes the error constant in 
Theorem \ref{thmmainstr} for $(T,\Sigma, \cF)=(T_pA, \Sigma (A), f)$.
It is proved by Banaszak, Gajda and Kraso\'n  
that $\mathrm{Im}\, \rho_{A,1}^{(p)}$ contains 
$\Sp_{2g} (\F_p)$ for sufficiently large $p$ 
if $A$ satisfies the following (i)--(iv).
\begin{enumerate}[{\rm (i)}]
\item  The abelian variety $A$ is  simple. 
\item  There is no endomorphism on $A$ 
defined over $\overline{K}$
except the multiplications by rational integers, namely 
$\End_{\overline{K}}(A) = \Z$.
\item  For any prime number $\ell$,  
the Zariski closure of the image of 
the $\ell$-adic representation 
\(
\rho^{(\ell)}_{A} \colon G_K \longrightarrow 
\Aut_{\Z_\ell} (V_\ell A) \simeq \GL_{2g} (\Q_\ell)
\)
is a connected algebraic group.
\item The dimension $g$ of $A$ is odd.
\end{enumerate}
See \cite{BGK} Theorem 6.16 for the cases when $\End_{\overline{K}}(A) = \Z$. 
(Note that in \cite{BGK}, they proved more general results. 
For details, see loc.\ cit..) 
\end{rem}

\begin{rem}\label{remabred}
Here, we shall describe the error terms $\nu_{v,n}$ 
in Theorem \ref{thmmainstr} for $T=T_pA$ 
in terminology related to the reduction of $A$ at $v$.  
Let  $\cA$ be the N\'eron model of $A$ over $\cO_K$. 
Take any finite place $v \in P_K$, and denote by $k_v$ 
the residue field of $\cO_{K_v}$. 
We put $A_{0,v}:=\cA \otimes_{\cO_K}k_v$, 
and define $A_{0,v}^0$ to be the identity component 
of $A_{0,v}$.  
Note that we have 
$H^0( K_v^{\unram}, A[p^\infty]) \simeq 
A_{0,v}(\overline{k}_v)$. (See \cite{ST} \S 1, Lemma 2.) 
By Chevalley decomposition, 
we have an exact sequence
\(
0 \longrightarrow T_w \times U_w 
\longrightarrow A_{0,v}^0
\longrightarrow B_v
\longrightarrow 0
\)
of group schemes over $k$,
where $T_w$ is a torus, 
$U_w$ is a unipotent group, 
and $B_v$ is an abelian variety. 
(For instance, see \cite{Co} Theorem 1.1 
and \cite{Wa} Theorem 9.5.)
In particular,   
if $U_w( \overline{k}_v)[p^\infty]=0$ 
(for instance, if $v$ does not lie above $p$), 
then the divisible part of 
$A_{0,v}( \overline{k}_v)[p^\infty]$
coincides with $A^0_{0,v}( \overline{k}_v)[p^\infty]$, 
and hence 
\begin{equation}\label{eqnured}
\nu_{v,n}=\ell_{\Z_p} 
(\pi_0 (A_{0,v})(k_v)[p^\infty]\otimes_\Z \Z/p^n \Z)
=\ell_{\Z_p} 
(\pi_0 (A_{0,v})(k_v)[p^n]), 
\end{equation}
where $\pi_0 (A_{0,v})$ denotes the group of 
the connected components of  $A_{0,v}$. 
(Note that the second equality holds since $\pi_0 (A_{0,v})$ 
is  finite.)
We can compute the error factors $\nu_{v,n}$ explicitly
if we know the structure of the reduction $A_{0,v}$ of $A$ 
at each finite place $v$. 
\end{rem}

\begin{ex}[Error factors for elliptic curves]\label{exellstr}
Now, we study the error factors $\nu_{v,n}$ in the setting of 
\cite{SY2} and \cite{Hi}.
We set $K=\Q$, and let 
$A$ be an elliptic curve with minimal 
discriminant $\Delta$.
Let $p$ be a prime number satisfying $\Full$
in Example \ref{exmnuimvanish}. 
We assume that $p$ is odd.  
Moreover, we also assume the following hypothesis. 
\begin{itemize}
\item If $p=3$, then $A$ does not have additive reduction at $p$. 
\item If $A$ has additive reduction at $p$, 
then $A(\Q_p)[p]=0$. 
\item If $A$ has split multiplicative reduction at $p$, 
$p$ does not divide $\ord_p (\Delta)$. 
\end{itemize}
(These hypotheses are assumed in
\cite{SY2} for $p>2$ and \cite{Hi}.)
Note that $\Sigma(A)$ is the set of places dividing 
$\infty p \Delta$. 
In this situation, we have 
$U_p( \overline{\F}_p)[p^\infty]=0$, 
where $U_p$ is the unipotent part of 
$A_{0,p}^0$.  
So, by (\ref{eqnured}), we obtain $\nu_{p, n}=0$ 
for each $n \in \Z_{\ge 0}$ since  by
\cite{Si} CHAPTER IV \S 9 Tate's algorithm 9.4,
\begin{itemize}
\item if $A$ has good reduction at $p$, then 
$A_{0,p}$ is connected; 
\item if $A$ has split multiplicative reduction,
then $\pi_0 (A_{0,p})(\F_p) \simeq \Z/\ord_p (\Delta) \Z$; 
\item if $A$ has non-split multiplicative reduction,
then $\pi_0 (A_{0,p})(\F_p) \simeq 0\ \text{or}\ \Z/2\Z$; 
\item if $A$ has additive reduction, 
then the order of $\pi_0 (A_{0,p})(\F_p)$ is prime to $p$. 
\end{itemize}

Let $\ell$ be a prime number distinct from $p$.
Then, we have $U_\ell ( \overline{\F}_\ell)[p^\infty]=0$. 
So, similarly to the above arguments, we obtain 
\begin{equation}\label{eqnuell}
\nu_{\ell, n} = 
\begin{cases}
\min \{ \ord_p (\ord_\ell (\Delta)), n) \} & 
\left(\begin{array}{l}
\text{if $A$ has split multiplicative} \\
\text{reduction at $\ell$}
\end{array}
\right), \\
0 & (\text{otherwise}).
\end{cases}
\end{equation}
Combining with Example \ref{exmnuimvanish}, 
Theorem \ref{thmmainstr} implies that
\[
\ord_p (h_n(A;p)) \ge 
2 \left( r_{\Z} (A)
- 1  \right) n
-\sum_{p \ne \ell\mid \Delta} \nu_{\ell, n},
\]
where $\nu_{\ell, n}$ is given by (\ref{eqnuell}).
This inequality coincides with that obtained in 
\cite{SY1}, \cite{SY2} and \cite{Hi} when $p$ is odd. 
In \cite{SY2}, Sairaiji and Yamauchi also treat the cases 
when $p=2$. 
Note that when $p=2$, an inequality 
following from Theorem \ref{thmmainstr} 
is weaker than that obtained in \cite{SY2}. 
Indeed, our $\nu_{p,n}$ is always non-negative by definition, 
but instead of it, in \cite{SY2}, they introduced 
a constant $\delta_2$ related to the local behavior of $A$ at $p=2$
which may become a negative integer.
\end{ex}

\subsection{RM and CM cases}\label{seccm}

Here we shall prove  Corollary \ref{corCM}. 
Let  
$K/K^+, p, \pi, A$ and $h_n (A;\pi)$ 
be as in Corollary \ref{corCM}. 
Take a subset $\Phi =\{\phi_1 , \dots , \phi_g \}
\subseteq \Gal (K/\Q)$ such that we have an isomorphism
\begin{equation}\label{eqstrCM}
\Lie (A/K) \simeq \bigoplus_{i=1}^g (K, \phi_i)
\end{equation} 
of modules over the ring 
\(
K \otimes_\Z\End (A)=
K \otimes_\Z \cO_K = \prod_{\sigma \in \Gal (K/\Q)} (K,\sigma)
\). 
Note that $\phi_1\vert_{K^+}, \dots, \phi_g\vert_{K^+}$
are distinct $g$ elements of $\Gal (K^+ /K)$.

Let us introduce  notation related to 
the formal group law.
Take any $\sigma \in \Gal (K/\Q)$, and 
denote by ${\sigma(\pi)}$  the place of $K$ 
corresponding to $\sigma(\pi) \cO_K$ (by abuse of notation). 
We put $k_{\sigma (\pi)}:=
\cO_{K}/\sigma (\pi) \cO_{K}=\F_p$. 
We define $\cA_{\sigma(\pi)}$ 
to be the N\'eron model of 
$A_{K_{\sigma (\pi)}}:=
A \otimes_K K_{\sigma (\pi)}$ over $\cO_{K_{\sigma(\pi)}}$, 
and $O_{\sigma (\pi), s}$ (resp.\ $O_{\sigma (\pi), \eta}$)  
the origin of the special (resp.\ generic) fiber of 
$\cA_{\sigma(\pi)}$.
For each $\star \in \{ s, \eta \}$, 
let $\mathfrak{m}_{\sigma(\pi),\star}$ 
be the maximal ideal of the local ring 
$\sO_{\widehat{\cA}_{\sigma(\pi)},\star}$. 
Note that $\sO_{{\cA}_{\sigma(\pi)},s}$ is 
a  regular local ring 
since $A$ has good reduction at $\sigma (\pi)$.
Let $s'_{0}=p, s'_{ 1}, \dots , s'_{g} 
\in \mathfrak{m}_{{\sigma(\pi)},s}$ be any 
regular system of parameters for the local ring 
$(\sO_{{\cA}_{\sigma(\pi)},O_{{\sigma(\pi)},s}}, 
\mathfrak{m}_{{\sigma(\pi)},s})$. 
Put $\mathfrak{n}:=\sO_{{\cA}_{{\sigma(\pi)},O_{{\sigma(\pi)},s}}} 
\cap \mathfrak{m}_{\sigma(\pi),\eta}$. 
Since we have the identity section 
$\Spec \cO_{K_{\sigma (\pi)}} \longrightarrow \cA_{\sigma (\pi)}$, 
it holds that
$\sO_{{\cA}_{\sigma(\pi)},O_{{\sigma(\pi)},s}}/\mathfrak{n}
\simeq \cO_{K_{\sigma (\pi)}}$. 
So, for each $i \in \Z$ with 
$1 \le i \le g$, we have a unique element 
$c_i \in \sigma (\pi) \cO_{K_{\sigma (\pi)}}$ 
such that $s_i:=s'_i-c_i \in \mathfrak{n}$.
Since $\cO_{K_{\sigma (\pi)}}$ is a DVR, 
and since we have 
\begin{align}
\mathfrak{n}/\mathfrak{n}^2 
\otimes_{\cO_{K_{\sigma (\pi)}}} k(\sigma (\pi))
& \simeq \coLie (A_{0,\sigma(\pi)}/k(\sigma (\pi))), \\
\mathfrak{n}/\mathfrak{n}^2 
\otimes_{\cO_{K_{\sigma (\pi)}}} K_{\sigma (\pi)}
& \simeq \coLie (A_{K_{\sigma (\pi)}}/K_{\sigma (\pi)}), 
\label{eqcoLie}
\end{align}
the  sequence $s_{ 1}, \dots , s_{g} $ forms 
a regular system of parameters for the regular local ring 
$(\sO_{{\cA}_{\sigma(\pi)},O_{{\sigma(\pi)},\eta}}, 
\mathfrak{m}_{{\sigma(\pi)}, \eta})$, and 
it holds that 
\begin{equation}\label{eqn/n2}
\mathfrak{n}/\mathfrak{n}^2 = \bigoplus_{i=1}^g  
\cO_{K_{\sigma (\pi)}} \bar{s}_i , 
\end{equation}
where $A_{0,\sigma(\pi)}$ denotes the special fiber 
of ${\cA}_{\sigma(\pi)}$, 
and $\bar{s}_i$ denotes the image of $s_i$ 
in $\mathfrak{n}/\mathfrak{n}^2$.

We denote by 
$\widehat{\cA}_{\sigma(\pi)}= \Spf 
\sO_{\widehat{\cA}_{\sigma(\pi)}}$ 
the completion of $\cA_{\sigma(\pi)}$ 
along $O_{\sigma(\pi),s}$, 
and by $\widehat{\mathfrak{m}}_{\sigma(\pi)}$
the maximal ideal of $\sO_{\widehat{\cA}_{\sigma(\pi)}}$.
Note that  
$s_{0}:=p, s_{ 1}, \dots , s_{g} $
forms a regular system of parameters for the 
complete regular local ring  
$(\sO_{\widehat{\cA}_{\sigma(\pi)}}, 
\widehat{\mathfrak{m}}_{\sigma(\pi)})$ 
We also note that $p$ splits completely  (unramified in particular) 
in $K/\Q$ by our assumption. 
So, we have  
$\sO_{\widehat{\cA}_{\sigma(\pi)}} 
= \cO_{K_{\sigma(\pi)}}[[s_1 , \dots , s_g]]$.
(See \cite{Ma} Theorem 29.7.)
For each $i \in \Z$ with $1 \le i \le g$, 
we define a formal power series 
$\sF_{A, \sigma(\pi),i}
\in \cO_{K_{\sigma(\pi)}}[[x_1, \dots , x_g, y_1, \dots, y_g]]$ 
by 
\[
\sF_{A, \sigma(\pi),i} := \mathrm{add}_{A, \sigma(\pi)}^{\sharp} (s_i), 
\]
where $\mathrm{add}_{A, \sigma(\pi)}^{\sharp} \colon 
\sO_{\widehat{\cA}_{\sigma(\pi)} }
\longrightarrow \sO_{\widehat{\cA}_{\sigma(\pi)}} 
\widehat{\otimes}_{\cO_{K_{\sigma(\pi)}}}
\sO_{\widehat{\cA}_{\sigma(\pi)}} 
=\cO_{K_{\sigma(\pi)}}[[x_1, \dots , x_g, y_1, \dots, y_g]]$
is the ring homomorphism corresponding to 
the group structure of the formal group scheme 
$\widehat{\cA}_{\sigma (\pi)}$. 
Note that since $s_1, \dots , s_g$ forms 
a regular system of parameters for the regular local ring 
$\sO_{{\cA}_{\sigma(\pi)},O_{{\sigma(\pi)},\eta}}$, 
the correction 
$\sF_{A, \sigma(\pi)}
=(\sF_{A, \sigma(\pi),i})_{i=1}^g$ 
is a $g$-dimensional commutative formal group law 
over $K_{\sigma(\pi)}$, 
and hence that over $\cO_{K_{\sigma(\pi)}}$. 
(For details,  see  Lemma C.2.1 in Appendix C of \cite{HS}.)
Let $\alpha \in \cO_K$ be any element, and    
$[\alpha]_{A, {\sigma(\pi)}}^{\sharp} \colon 
\sO_{\widehat{\cA}_{\sigma(\pi)}} 
\longrightarrow \sO_{\widehat{\cA}_{\sigma(\pi)}}$ 
the ring homomorphism corresponding to  
the multiplication-by-$\alpha$ endomorphism on 
the formal group scheme $\widehat{\cA}_{\sigma (\pi)}$.  
For each $i \in \Z$ with $1 \le i \le g$,  
we put 
\[
[\alpha]_{A, {\sigma(\pi)},i}(s_1 , \dots , s_g)
:= [\alpha]_{A, {\sigma(\pi)}}^\sharp (s_i) \in 
\sO_{\widehat{\cA}_{\sigma(\pi)}} 
= \cO_{K_{\sigma(\pi)}}[[s_1 , \dots , s_g]]. 
\]

\begin{lem}
\label{lemlocparam}
There exists a regular system of parameters 
$s_{ 0}=p, s_{ 1}, \dots , s_{ g} 
\in \mathfrak{m}_{\sigma(\pi),s}$
such that for any $i \in \Z$ with $1 \le i \le g$ 
and any $\alpha \in \cO_K$,  it holds that $s_i \in \mathfrak{n}$, 
and \[
[\alpha]_{A, {\sigma(\pi)},i}(s_1 , \dots , s_g)
\equiv \phi_i (\alpha) s_i \mod \mathfrak{n}^2.
\] 
\end{lem}

\begin{proof}
Since we have a direct product decomposition 
\[
\cO_{K_{\sigma (\pi)}} \otimes_{\Z} \cO_K= \Z_p  
\otimes_{\Z} \cO_K= 
\varprojlim_{n} \cO_K/p^n \cO_K
=\prod_{\tau \in \Gal (K/\Q)} 
\cO_{K_{\tau (\pi)}}
\] 
by  Chinese remainder theorem, 
and by (\ref{eqstrCM}),   
(\ref{eqcoLie}) and (\ref{eqn/n2}), 
we can take an $\cO_{K_{\sigma (\pi)}}$-basis 
$\bar{s}_1, \cdots,  \bar{s}_g$ of $\mathfrak{n}/\mathfrak{n}^2$ 
such that for each $i \in \Z$ with $1 \le i \le g$ and 
each $\alpha \in \cO_K$, the element 
$\bar{s}_i$ is an eigenvector of 
the multiplication-by-$\alpha$ map $[\alpha]$
attached to the eigenvalue $\phi_i (\alpha)$. 
We take any lift $s_1, \cdots,  s_g \in \mathfrak{n}$ of 
$\bar{s}_1, \cdots,  \bar{s}_g$.
Then the sequence
$s_{ 0}=p, s_{ 1}, \dots , s_{ g}$ 
is the one as desired. 
\end{proof}

From now on, let the parameters $s_1, \dots , s_g$ 
be as in \ref{lemlocparam}. 
We define 
\[
\sF_{A, {\sigma(\pi)}} 
( \sigma (\pi) \cO_{K_{\sigma (\pi)}})=
\left(
(\sigma (\pi) \cO_{K_{\sigma (\pi)}})^g, 
\sF_{A, {\sigma(\pi)}} \right)\]
to be the set $(\sigma (\pi) \cO_{K_{\sigma (\pi)}})^g$ equipped with 
a group structure defined by the formal group law $\sF_{A, {\sigma(\pi)}}$. 
Note that by the scalar multiplication defined by 
the collection $([\alpha]_{A, {\sigma(\pi)},i})_{i=1}^g$ 
of the power series,   
we can regard  $\sF_{A, {\sigma(\pi)}} 
( \sigma( \pi ) \cO_{K_{\sigma (\pi)}})$
as  an $\cO_K$-module.
By Lemma \ref{lemlocparam}, the following holds.

\begin{cor}\label{corformal}
Let  $\sigma \in \Gal (K/\Q)$ any element.
Then, we have 
\[
\sF_{A, {\sigma(\pi)}} 
\left(
\sigma(\pi) \cO_{K_{\sigma (\pi)}}
\right)
\otimes_{\cO_{K}} \varinjlim_{n>0} 
\pi^{-n} \cO_{K}/\cO_{K}
\simeq 
\begin{cases}
\Q_p/\Z_p & (\text{if $\sigma \in \Phi$}), \\
0 & (\text{if $\sigma \notin \Phi$}).
\end{cases}
\]
\end{cor}

Let $A_{0,\sigma(\pi)}$ be the special fiber 
of ${\cA}_{\sigma(\pi)}$.
Since $A$ has good reduction at $\sigma (\pi)$, 
we can define 
the reduction map 
\(
\mathrm{red}_{A,\sigma(\pi)}
\colon A(K_{\sigma (\pi)}) 
\longrightarrow A_{0,v}(k_{\sigma (\pi)}) 
\).

\begin{lem}[For instance, see \cite{HS} Theorem C.2.6]\label{lemredformal}
For any finite place $v \in P_K$, 
the $\cO_K$-module 
$\Ker \mathrm{red}_{A,\sigma(\pi)}$
is isomorphic to 
$\sF_{A, {\sigma(\pi)}} (\pi_v \cO_{K_{\sigma (\pi)}})$. 
\end{lem}

\begin{rem}
Note that 
\cite{HS} Theorem C.2.6 says that 
$\Ker \mathrm{red}_{A,\sigma(\pi)}$
is isomorphic to 
$\sF_A (\pi_v \cO_{K_v})$
only as a group, but  it is easy to verify that 
the group isomorphism constructed in the proof of 
\cite{HS} Theorem C.2.6 preserves the scalar action of $\cO_K$. 
\end{rem}

We define $T_\pi A:=\varprojlim_n A[\pi^n]$. 
Let $\Sigma(A)$ be a subset of $P_K$ consisting of 
all places dividing $p \infty$ and all places 
where $A$ has bad reduction.
Since $T_\pi A$ is regarded as a $\Z_p$-submodule of $T_p A$, 
the action of $\Gal (\overline{K}/K)$ on $T_\pi A$ 
is unramified outside $\Sigma(A)$. 
Since $p$ splits completely in $K$, 
the $\Z_p$-module 
$T_\pi A$ is free  of rank 
$2 \dim/[K:\Q]=2/[ K : K^+ ]$. 
Note that $T_\pi A$ satisfies $\Abs$ and $\NT$
by our assumption.  

Take a finite place  
$v \in P_K$ above a prime number $\ell$.
Since $H^1(K_v,A[\pi^\infty])$ is a direct summand of 
$H^1(K_v,A[p^\infty])$ consisting of elements 
annihilated by $\pi^n$ for some $n \in \Z_{\ge 0}$, 
Proposition \ref{corclvsfin} implies that  
\begin{align*}
H^1_f(K_v,A[\pi^\infty])& =H^1_f(K_v,A[p^\infty]) 
\cap H^1(K_v,A[\pi^\infty])
=H^1_f(L_v,A[p^\infty]) [\pi^\infty] \\
& \simeq 
(A(K_v) \otimes_{\Z_p}\Q_p/\Z_p)[\pi^\infty] \simeq 
A(K_v) \otimes_{\cO_{K}} K_\pi/
\cO_{K_\pi}.
\end{align*}
By isomorphisms (\ref{eqstrthm}) (for $v \nmid p$), 
Corollary \ref{corformal} and
Lemma \ref{lemredformal}, the following holds. 

\begin{prop}\label{propfinclCM}
Let $v \in P_K$ be a finite place. 
Then, we have 
\[
H^1_f(K_v,A[\pi^\infty])=
\begin{cases}
\Q_p/\Z_p & (\text{if $v= \sigma(\pi)$ 
for some $\sigma \in \Phi$}), \\
0 & (\text{otherwise}).
\end{cases}
\]
\end{prop}

\begin{proof}[Proof of Corollary \ref{corCM}]
Note that the Selmer group
$\Sel_f(K,A[\pi^\infty])$ is 
a direct summand of the $\Z_p$-module
$\Sel_f(K,A[p^\infty])$ consisting of 
elements annihilated by $\pi^n$
for some $n \in \Z_{\ge 0}$. 
So, we have
\(
\Sel_f(K,A[\pi^\infty]) \simeq 
\Sel_f(K,A[p^\infty]) \otimes_{\cO_K} K_\pi/
\cO_{K_\pi}
\).
Combining with the short exact sequence (\ref{eqexseqSel}), 
this implies that $\Sel_f(K,A[\pi^\infty])$ has 
a $\Z_p$-submodule isomorphic to 
$A(K) \otimes_{\cO_K} 
K_\pi/
\cO_{K_\pi}$. 
Thus it holds that  
\[
\corank_{\Z_p} 
\Sel_f(K,A[\pi^\infty]) 
\ge r_{\cO_K} (A) := \dim_K (A(K) \otimes_{\cO_K} K).
\] 
Hence by Theorem \ref{thmmain} for 
$(T,\Sigma, \cF)=(T_\pi A, \Sigma (A), f)$ and 
Proposition \ref{propfinclCM}, we obtain 
the assertion of Corollary \ref{corCM}, that is,  
\(
h_{n}(A;\pi) \succ 2 [K : K^+]^{-1} 
\left(
r_{\cO_K} (A)  -g
\right) n 
\). 
\end{proof}

\begin{rem}
Here, let $p$ be an odd prime number, 
and consider the image of $G_{1,0};=\Gal (K_1/K)$ 
in $\Aut_{\F_p}(A[\pi])$. 
First, let  $K$ be a CM field.  
By Banaszak--Gajda--Kraso\'n's work 
(\cite{BGK} Theorem 6.16 for the cases 
then $\End_{\overline{K}} (A)=\cO_K$) 
it is known that  
the image of $G_{1,0}:=\Gal (K_1/K)$ 
in $\Aut_{\F_p}(A[\pi])= \GL_2 (\F_p)$ 
contains $\SL_2 (\F_p)$ if $A$ is principally polarized, and if 
$A$ satisfies (i) and (iii) in Remark \ref{remnuimforab}. 
So, if $A$ satisfies the conditions 
(i) and (iii) in Remark \ref{remnuimforab}, 
then $T_\pi A$ satisfies $\Abs$ and $\NT$, 
and we can take $\nuim=0$ by Remark \ref{remnuimvanish}. 

Next, let  $K$ be a CM field.  
Then  $T_\pi A$ obviously 
satisfies the condition $\Abs$ 
since $\dim_{\F_p} A[\pi]=1$. 
Moreover, if  $K$ is a CM field, 
then $G_{1,0}=\Gal (K_1/K) \simeq \F_p^\times$. 
(See, for instance \cite{Ro} Proposition 3.1.)
So, in this situation, the $\pi$-adic Tate module  
$T_\pi A$ also satisfies the condition $\NT$, 
and we can take $\nuim =0$ by Remark \ref{remnuimvanish}.  
\end{rem}

\end{document}